\documentclass{article}
\usepackage{amssymb}
\usepackage{amsmath}
\usepackage{amsfonts}
\usepackage{graphicx}

\newtheorem{theorem}{Theorem}[section]
\newtheorem{lemma}[theorem]{Lemma}

\newtheorem{remark}[theorem]{Remark}

\newcommand{\qed}{\hfill \rule{2.3mm}{2.3mm}}

\newcommand{\al}{\alpha}

\newcommand{\la}{\lambda}

\newcommand{\ve}{\varepsilon}
\newcommand{\vp}{\varphi}

\newcommand*{\esssup}{\operatorname*{ess\phantom{|}\!sup}}
\newcommand*{\essinf}{\operatorname*{ess\phantom{|}\!inf}}

\newcommand{\bee}{\begin{equation}}
\newcommand{\ee}{\end{equation}}

\newcommand{\R}{\mathbb{R}}

\newcommand{\M}{\mathbb{M}}
\newcommand{\MM}{{\cal M}_r}

\newcommand{\x}{x/\varepsilon}
\newcommand{\y}{y/\varepsilon}

\newcommand{\A}{{\cal A}}

\newcommand{\reff}[1]{(\ref{#1})}   
\begin{document}

\title{\bf 
Nonlinear non-periodic homogenization: 
Existence, local uniqueness and estimates 
} 
\newcounter{thesame}
\setcounter{thesame}{1}
\author{
	Lutz Recke \thanks{Humboldt University of Berlin, Institute of Mathematics, Rudower Chausee 25, 12489 Berlin, Germany.
		{\small   E-mail:
			{\tt lutz.recke@hu-berlin.de}}
}}

\date{}

\maketitle

\begin{abstract}
\noindent
We consider  periodic homogenization with localized defects of
boundary value problems for semilinear 
ODE systems of the type
$$
\Big((A(\x)+B(\x))u'(x)+c(x,u(x))\Big)'= d(x,u(x)) \mbox{ for } x \in (0,1),\; u(0)=u(1)=0.
$$
Here $A \in L^\infty(\R;\M_n)$ is 1-periodic, $B \in L^\infty(\R;\M_n)\cap
L^1(\R;\M_n)$, $A(y)$ and $A(y)+B(y)$ are positive definite uniformly with respect to~$y$, $c(x,\cdot),d(x,\cdot)\in C^1(\R^n;\R^n)$ and $c(\cdot,u),d(\cdot,u) \in L^\infty((0,1);\R^n)$.
For small $\ve>0$ we show existence of weak solutions $u=u_\ve$ as well as their local uniqueness for $\|u-u_0\|_\infty \approx 0$,
where $u=u_0$ is a given
non-degenerate weak
solution to the homogenized problem
$$
\left(\left(\int_0^1A(y)^{-1}dy\right)^{-1}
u'(x)+c(x,u(x))\right)'= d(x,u(x)) \mbox{ for } x \in (0,1),\; u(0)=u(1)=0.
$$
Further, we show that $\|u_\ve-u_0\|_\infty\to 0$
and, if $c(\cdot,u_0(\cdot))
\in W^{1,\infty}((0,1);\R^n)$, that  $\|u_\ve-u_0\|_\infty=O(\ve)$
for $\ve \to 0$.
Moreover, all these statements are true, roughly speaking, uniformly with respect to the localized defects $B$.

The main tool of the proofs is an abstract result of implicit function theorem type
which has been tailored for applications to nonlinear singular perturbation and homogenization problems.
\end{abstract}



\section{Introduction}
\setcounter{equation}{0}
\setcounter{theorem}{0}

This paper concerns periodic homogenization with localized defects in the sense of \cite{Blanc}. More exactly, we consider Dirichlet problems for semilinear second-order ODE systems in divergence form of the type
\bee
\label{ODE}
\left.
\begin{array}{l}
\Big((A(\x)+B(\x))u'(x)+c(x,u(x))\Big)'= d(x,u(x)) \mbox{ for } x \in (0,1),\\
u(0)=u(1)=0.
\end{array}
\right\}
\ee
For small $\ve>0$ we look for solutions $u:[0,1]\to \R^n$ to \reff{ODE}, and we suppose that 
\begin{eqnarray}
\label{aass}
&&A \in L^\infty(\R;\M_n) \mbox{ is 1-periodic,} \\
\label{mon}
&&\essinf\{A(y)u\cdot u:\; y \in \R, u \in \R^n, \|u\|=1\}>0,\\
\label{cass}
&&u \in \R^n \mapsto (c(\cdot,u),d(\cdot,u)) \in  L^\infty((0,1);\R^n)^2 \mbox{ is $C^1$-smooth.}
\end{eqnarray}
Here $\M_n$ is the space of all real $n\times n$-matrices, $(u,v)\in \R^n\times\R^n \mapsto u\cdot v \in \R$
is the Euclidean scalar product in 
$\R^n$, and
$$
\|u\|:=\sqrt{u\cdot u} \mbox{ for } u \in \R^n, \mbox{ and }
\|M\|:=\max\{\|Mu\|:\, u \in \R^n,\, \|u\|\le 1\}
$$
are the Euclidean norms in $\R^n$ and $\M_n$, respectively.
Further,
$$
\|u\|_\infty:=\esssup\{\|u(x)\|:\; x \in (0,1)\} 
$$
is the norm in the Lebesgue space $L^\infty((0,1);\R^n)$(which is used, for example, in \reff{cass}), $\|u\|_\infty+\|u'\|_\infty$ is the norm in the Sobolev space $W^{1,\infty}((0,1);\R^n)$(which is used in Theorem \ref{main}(ii) below), and
$$
\|M\|_\infty:=\esssup\{\|M(x)\|:\; x \in \R\} 
\mbox{ and }
\|M\|_1:=\int_{-\infty}^\infty\|M(x)\|dx
$$
are the norms in $L^\infty((0,1);\M_n)$ and $L^1((0,1);\M_n)$, respectively.
Finally, for $r>1$ we define
$$
{\cal B}_r:=\left\{B \in 
L^\infty(\R;\M_n)\cap L^1(\R;\M_n):\, \|M\|_\infty+\|M\|_1\le r\right\}
$$
and
\bee
\label{Mdef}
{\cal M}_r:=\left\{B \in {\cal B}_r:
\,
(A(y)+B(y))u\cdot u \ge \frac{\|u\|}{r}
\mbox{ for almost all } y \in \R \mbox{ and all } u \in \R^n\right\}.
\ee

We are going to prove existence and local uniqueness of
weak solutions $u$ to \reff{ODE} with $\ve \approx 0$, $B \in \MM$ and $\|u-u_0\|_\infty\approx 0$, where $u_0$ is a given non-degenerate weak solution to the homogenized boundary value problem
\bee
\label{hombvp}
\left.
\begin{array}{l}
\Big(A_0u'(x)+c(x,u(x)\Big)'= d(x,u(x)) \mbox{ for } x \in (0,1),\\
u(0)=u(1)=0
\end{array}
\right\}
\ee
with
\bee
\label{Anulldef}
A_0:=\left(\int_0^1A(y)^{-1}dy\right)^{-1}.
\ee
Remark that the assumption \reff{aass}
implies that the matrices $A(y)$ 
and, for $B\in \MM$, the matrices $A(y)+B(y)$ 
are positive definite uniformly with respect to $y \in \R$. Hence, these matrices are invertible, and their inverses also
are positive definite uniformly 
 with respect to $y \in \R$, and
\bee
\label{inv}
\esssup_{y \in \R}\|A(y)^{-1}\|
+\esssup_{y \in \R}\sup_{B\in \MM}
\|(A(y)+B(y))^{-1}\|
<\infty.
\ee
Also the matrix $\int_0^1A(y)^{-1}dy$ 
is positive definite and, hence, invertible. In particular,
the definition \reff{Anulldef} is correct. 

As usual, a vector function $u\in W_0^{1,2}((0,1);\R^n)$ is called weak solution to \reff{ODE} if 
for all $\vp \in W_0^{1,2}((0,1);\R^n)$
it
satisfies 
the variational equation
\bee
\label{vareq}
\displaystyle\int_0^1\Big(\Big((A(\x)+B(\x))u'(x)+c(x,u(x))\Big)\cdot\vp'(x)
+d(x,u(x))\cdot\vp(x)\Big)dx=0,
\ee
and similar for the homogenized boundary value problem \reff{hombvp} and for its linearization in $u=u_0$, which is
\bee
\label{linhombvp}
\left.
\begin{array}{l}
\Big(A_0u'(x)+\partial_uc(x,u_0(x))u(x)
\Big)'= \partial_ud(x,u_0(x))u(x)
 \mbox{ for } x \in (0,1),\\
u(0)=u(1)=0.
\end{array}
\right\}
\ee
Here $\partial_uc(x,u_0(x))\in \M_n$
is the Jacobi matrix of the vector function $u\in \R^n \mapsto  c(x,u)\in \R^n$ in the point $u=u_0(x)$, and similarly 
$\partial_ud(x,u_0(x))\in \M_n$.

Now we formulate our main result:
\begin{theorem} 
\label{main}
Suppose \reff{aass}-\reff{cass}, and
let $u=u_0$ be a weak solution to \reff{hombvp} such that
\reff{linhombvp}
does not have weak solutions $u\not=0$.

Then for any $r>1$ there exist
$\ve_r>0$ and  $\delta_r>0$
such that for all $\ve \in (0,\ve_r]$
and $B \in {\cal M}_r$
the following is true:

(i) There exists exactly one 
weak solution $u=u_{\ve,B}$ 
to \reff{ODE} with $\|u-u_0\|_\infty \le \delta$. Moreover,
$\|u_{\ve,B}-u_0\|_\infty\to 0$ for  $\ve \to 0$
uniformly with respect to $B \in {\cal M}_r$.

(ii) If $c(\cdot,u_0(\cdot))\in W^{1,\infty}((0,1);\R^n)$,
then $\|u_{\ve,B}-u_0\|_\infty=O(\ve)$ for  $\ve \to 0$ uniformly with respect to
$B \in {\cal M}_r$
\end{theorem}

\begin{remark}
If for $x \in (0,1)$ and $u \in \R^n$ we have $c(x,u)=\al_1(x)c_1(u)+\ldots+\al_m(x)c_m(u)$ with $\al_l \in L^\infty(0,1)$ and $c_l\in C^1(\R^n;\R^n)$, then this function $c$ satisfies assumption \reff{cass} of Theorem~\ref{main}. Moreover, if $\al_l \in W^{1,\infty}(0,1)$, then the assumption  $c(\cdot,u_0(\cdot))\in W^{1,\infty}((0,1);\R^n)$
of Theorem~\ref{main}(ii) is satisfied also (because $u_0 \in W^{1,\infty}((0,1);\R^n)$, cf. Lemma \ref{prep} below).
\end{remark}

\begin{remark}
In \cite{I} is proven a result of the type of Theorem \ref{main} for quasilinear ODE systems of the type
\bee
\label{ODEI}
a(x,\x,u(x),u'(x))'=b(x,\x,u(x),u'(x))
\mbox{ for } x \in (0,1).
\ee
On the one hand, this ODE system is much more general then that in \reff{ODE}. But on the other hand, in \reff{ODEI} it is supposed that the vector functions $a(x,\cdot,u(x),u'(x))$ and $b(x,\cdot,u(x),u'(x))$ are periodic, i.e. 
\cite{I} concerns periodic homogenization, while we do not suppose that the matrix function $B$ is periodic, i.e. Theorem \ref{main} concerns non-periodic homogenization
(or, more exactly, periodic homogenization with localized defects, cf. \cite[Section 2.3]{Blanc}).
Moreover, in \reff{ODEI} it is supposed that the functions $a(x,\cdot,u(x),u'(x))$ and $b(x,\cdot,u(x),u'(x))$ are continuous, while we do not suppose that the matrix functions $A$ and $B$ and the vector functions $c(\cdot,u)$ and $d(\cdot,u)$ are continuous.
\end{remark}

\begin{remark}
In \cite{II} is proven a result of the type of Theorem \ref{main} for periodic homogenization of semilinear elliptic PDE systems of the type
\bee
\label{PDES}
\nabla \cdot \Big(A(\x)\nabla u(x)+c(x,u(x))\Big)=d(x,u(x))\mbox{ for } x \in \Omega.
\ee
On the one hand, there are no localized defects in \reff{PDES}. On the other hand, the proofs in \cite{II} are more difficult because we have to work with  $\ve$-depending approximate solutions,
i.e. functions, which satisfy \reff{PDES} approximately for $\ve \to 0$, 
which are close to $u_0$ for $\ve \to 0$, and
which are constructed by using so-called correctors.
In the present paper (as well as in \cite{I}) there is no need for using correctors, because the boundary value problem \reff{ODE} can be transformed into the 
integral equation
\reff{ueq}, and $u_0$ itself is an 
approximate solution to \reff{ueq}.
Accordingly, in \cite{II} there is used a version of Theorem \ref{ift}
where $u_0$ is replaced by a family of $\ve$-depending approximate solutions.
\end{remark}

\begin{remark}
\label{sufficient}
The assumption of Theorem \ref{main}, that there do not exist weak solutions  $u\not=0$ to \reff{linhombvp}, is rather implicit. But there exist simple explicit sufficient conditions for it. For example, if
not only the matrix $A_0$ is positive definite, 
(this follows from
\reff{mon}
and \reff{Anulldef}), 
but also the matrices $\partial_{u}d(x,u_0(x))$, 
and if $\|\partial_{u}c(\cdot,u_0(\cdot))\|_\infty$ is sufficiently small, then there do not exist nontrivial weak solutions to \reff{linhombvp}.
\end{remark}

\begin{remark}
\label{Neuk}
$L^\infty$-estimates of the homogenization error $u_\ve-u_0$ exist, to the best of our knowledge, for linear periodic homogenization problems with smooth coefficients only: For  scalar ODEs
of the type
$
\big(a(\x)u'(x)\big)'=f(x) 
$
in \cite[Section 1]{N}, for scalar ODEs with stratified structure of the type
$
\big(a(x,\rho(x)/\ve)u'(x)\big)'=f(x)
$
in \cite[Theorem 1.2]{YXu}.
For $L^\infty$ periodic homogenization error estimates for scalar linear elliptic PDEs of the type
$
\mbox{\rm div}\, a(\x) \nabla u(x)=f(x)
$
see, e.g. \cite[Chapter 2.4]{Ben} and \cite{He} and for linear elliptic systems \cite[Theorem 7.5.1]{Shen}.

For periodic homogenization of linear equations with $L^\infty$-coefficients see \cite[Theorems 6.1 and 6.3]{Ben} and \cite[Theorem 1.2]{YXu}.

What concerns existence and local uniqueness
for nonlinear periodic homogenization problems
(without assumption of global uniqueness), besides \cite{I} and \cite{II} we know only the result 
\cite{Bun}  for scalar semilinear elliptic PDEs of the type
$
\mbox{\rm div}\, a(\x) \nabla u(x)=f(x)g(u(x)),
$
where the nonlinearity $g$ is supposed to have a sufficiently small local Lipschitz constant (on an appropriate bounded interval). Let us mention also \cite{Lanza1,Lanza2}, where existence and local uniqueness for a periodic homogenization problem for the linear Poisson equation with nonlinear Robin boundary conditions is shown. There the specific structure of the problem (no highly oscillating diffusion coefficients) allows to apply the classical implicit function theorem.
\end{remark}

~\\ 

Our paper is organized as follows: 
In Section \ref{secabstract} we consider abstract nonlinear parameter depending equations of the type
\bee
\label{intrabstract}
u=F_\ve(u).
\ee
Here $\ve\ge 0$ is the parameter. We prove a result on existence and local uniqueness of a family of solutions $u=u_\ve \approx u_0$ to \reff{intrabstract} with $\ve \approx 0$, where $u_0$ is an approximate solution to \reff{intrabstract}, i.e. an element with
$F_\ve(u_0)\to u_0$ for $\ve \to 0$, and we estimate the norm of the error $u_\ve-u_0$ by the norm of the discrepancy $u_0-F_\ve(u_0)$.
Moreover, we show that all these results are true uniformly with respect to a second parameter (which is the localized defect $B$ in the application to problem \reff{ODE}), if the equation
\reff{intrabstract} depends appropriately on this second parameter.

This type of generalized implicit function theorems has been  applied to singularly perturbed nonlinear ODEs and to elliptic and parabolic PDEs in \cite{Butetc,But2022,Fiedler,
Magnus1,Magnus2,NURS,
OmelchenkoRecke2015,Recke2022,
ReckeOmelchenko2008}) as well as to periodic homogenization of nonlinear ODEs (in \cite{I}) and elliptic PDEs (in \cite{II}). 
The proofs of all these results are based on the generalized implicit function theorem
of R.J. Magnus \cite[Theorem 1.2]{Magnus1} and on several of its modifications (see, e.g.
Theorem \ref{ift} below).
Contrary to the classical implicit function theorem it is not supposed that for all $u \approx u_0$
the linearized operators $F'_\ve(u)$ converge for $\ve \to 0$ with respect to the uniform operator norm. And, indeed, in the applications to singularly perturbed problems as well as to  homogenization problems they do not converge for $\ve \to 0$ with respect to the uniform operator norm (cf. Remark \ref{uniform} below).
Remark that in the classical implicit function theorem one cannot omit, in general, the assumption, that
for all $u \approx u_0$ the linearized operators
$F_\ve'(u)$ converges for $\ve \to 0$ with respect to the uniform operator norm (cf. \cite[Section 3.6]{Katz}).


In Section \ref{sec3} we prove Theorem \ref{main} by means of the results of Section \ref{secabstract}. For that reason we transform the boundary value problem 
\reff{ODE} into the 
integral equation
\reff{ueq}, and for that integral equation
we introduce in \reff{abdef} an abstract setting
of the type
\reff{intrabstract}.
For that abstract setting we have to verify the assumptions \reff{compass}-
\reff{linlimass} of Theorem \ref{ift}, and we do this in the 
Subsections~\ref{subsec:Feps}-\ref{subsec:infass1}.

\section{An abstract result of implicit function theorem type}
\label{secabstract}
\setcounter{equation}{0}
\setcounter{theorem}{0}
In this section we formulate and prove Theorem \ref{ift} below. 
\begin{theorem}
\label{ift}
Let be given Banach spaces $U$ and $V$, a set $\Lambda$, a family of differentiable maps $F_{\ve,\la}:U \to V$ (with $\ve> 0$ and $\la \in \Lambda$ as family parameters), a differentiable map $F_0:U \to V$
 and an element $u_0 \in U$
such that the following is true:
\begin{eqnarray}
\label{compass}
&&\mbox{$V$ is compactly embedded into $U$,}\\
\label{Fnullass}
&&u_0=F_0(u_0),\\
\label{linass}
&&\ker(I-F'_0(u_0))=\{0\},\\
\label{supass}
&&\sup\{\|F'_{\ve,\la}(u_0)u\|_V:\; \ve >0,\, \la \in \Lambda,
\; u \in U, \; \|u\|_U\le 1\}<\infty,\\
\label{linlimass1}
&&\lim_{\ve+\|u\|_U\to 0} \sup_{\|v\|_U\le 1}\|(F'_{\ve,\la}(u_0+u)-F'_\ve(u_0))v\|_U=0 \mbox{ uniformly with respect to } \la \in \Lambda,\\
\label{limass}
&&\lim_{\ve\to 0} \|F_{\ve,\la}(u_0)-F_0(u_0)\|_U=0 \mbox{ uniformly with respect to } \la \in \Lambda,\\
\label{linlimass}
&&\lim_{\ve\to 0} \|(F'_{\ve,\la}(u_0)-F'_0(u_0))u\|_U=0
\mbox{ for all } u \in U \mbox{ uniformly with respect to } \la \in \Lambda.
\end{eqnarray}

Then there exist $\ve_0>0$ and $\delta >0$ such that for all $\ve \in (0,\ve_0]$ and $\la \in \Lambda$
there exists exactly one  $u = u_{\ve,\la} \in U$ with $u=F_{\ve,\la}(u)$ and $\|u-u_0\|_U \le\delta$. Moreover, 
there exists $\rho>0$ such that
\bee
\label{abest}
\|u_{\ve,\la}-u_0\|_U \le \rho\|F_{\ve,\la}(u_0)-F_0(u_0)\|_U \mbox{ for all }
\ve \in (0,\ve_0] \mbox{ and } \la \in \Lambda.
\ee
\end{theorem}
{\bf Proof } In a first step we show that there exists $\ve_0>0$ such that
\bee
\label{infass}
\inf
\left\{\|(I-F'_{\ve,\la}(u_0))u\|_U:\,
\ve \in (0,\ve_0],\, \la \in \Lambda,\,u \in U,\,\|u\|_U=1\right\}
=: \al>0.
\ee
Suppose the contrary. Then there exist sequences $\ve_1,\ve_2,\ldots>0$,
$\la_1,\la_2,\ldots \in \Lambda$ and
$u_1,u_2,\ldots \in U$ such that
\bee
\label{convass}
\lim_{l \to \infty}\Big(\ve_l+\|(I-F'_{\ve_l\la_l,}(u_0))u_l\|_U\Big)=0,
\ee
but
\bee
\label{eqass}
\|u_l\|_U=1 \mbox{ for all } l.
\ee
Because of \reff{supass} and  \reff{eqass} the sequence $F'_{\ve_1,\la_1}(u_0)u_1,F'_{\ve_2,\la_2}(u_0)u_2,\ldots$ is bounded in $V$. Hence, \reff{compass} implies that without loss of generality we may assume that there exists $u_* \in U$ such that 
\bee
\label{ustar}
\lim_{l \to \infty}\|u_*-F'_{\ve_l,\la_l}(u_0)u_l\|_U=0.
\ee
 Therefore \reff{convass} yields that
\bee
\label{starconv}
\lim_{l \to \infty}\|u_l-u_*\|_U=0,
\ee
and \reff{supass}, \reff{ustar} and 
\reff{starconv} imply 
$\|u_*-F'_{\ve_l,\la_l}(u_0)u_*\|_U\to 0$
for $l \to \infty$.
Hence \reff{linlimass} implies that
$u_*=F'_{0}(u_0)u_*$, and \reff{linass} yields that $u_*=0$. But this contradicts  \reff{eqass} and \reff{starconv}.

Now, in a second step we solve the equation $u=F_{\ve,\la}(u)$ for $\ve \approx 0$, $\la \in \Lambda$ and $\|u-u_0\|_U\approx 0$.

Because the maps $F_{\ve,\la}$ are differentiable from $U$ into $V$, their linearizations $F'_{\ve,\la}(u_0)$ are bounded from $U$ into $V$ and, because of \reff{compass}, compact from $U$ into $U$. 
Therefore \reff{infass} implies that 
for $\ve \in (0,\ve_0]$ and $\la \in \Lambda$ the linear operators $I-F'_{\ve,\la}(u_0)$ are isomorphisms from $U$ onto $U$, and 
\bee
\label{invest}
\|(I-F'_{\ve,\la}(u_0))^{-1}u\|_U\le \frac{1}{\al}\|u\|_U
\mbox{ for all } \ve \in (0,\ve_0], \la \in \Lambda
\mbox{ and } u \in U.
\ee

Obviously, for $\ve \in (0,\ve_0]$, 
$\la \in \Lambda$
and $u \in U$ we have $u=F_{\ve,\la}(u)$ if and only if $u$ is a solution
\bee
\label{fp}
u=G_{\ve,\la}(u)
\ee
with
$$
G_{\ve,\la}(u):=u-(I-F'_{\ve,\la}(u_0))^{-1}(u-F_{\ve,\la}(u))=(I-F'_{\ve,\la}(u_0))^{-1}(F_{\ve,\la}(u)-F'_{\ve,\la}(u_0)u).
$$
We are going to solve \reff{fp} by means of Banach's fixed point theorem for $\ve \approx 0$, $\la \in \Lambda$ and $\|u-u_0\|_U\approx 0$.
For $r>0$ denote 
$ 
\mathbb{B}_r
:=\{u \in U:\; \|u-u_0\|_U \le r\}.
$
We have to show that for sufficiently small $\ve>0$ and $r>0$ the map $G_{\ve,\la}$ is strictly contractive from the closed ball $\mathbb{B}_r$
into itself.

In order to verify the strict contractivity of $G_{\ve,\la}$ we take 
$\ve \in (0,\ve_0]$, $\la \in \Lambda$ and $u_1,u_2 \in U$ and estimate as follows:
\begin{eqnarray*}
\lefteqn{
\|G_{\ve,\la}(u_1) - G_{\ve,\la}(u_2)\|_U = \left\|(I-F'_{\ve,\la}(u_0))^{-1}\Big(
F_{\ve,\la}(u_1)-F_{\ve,\la}(u_2)-F'_{\ve,\la}(u_0)(u_1-u_2)\Big)\right\|_U}\nonumber\\
&&=\left\|(I-F_{\ve,\la}'(u_0))^{-1}\int_0^1\Big(F'_{\ve,\la}(s u_1 + (1-s) u_2)-F'_{\ve,\la}(u_0)\Big)ds\, (u_1 - u_2)\right\|_U\nonumber\\
&&\le \frac{1}{\al}\int_0^1
\left\|\Big(F'_{\ve,\la}(s u_1 + (1-s) u_2)-F'_{\ve,\la}(u_0)\Big)(u_1 - u_2)\right\|_Uds.
\end{eqnarray*}
Here we used \reff{invest}.
Because of assumption \reff{linlimass1} there exist $\ve_1 \in (0,\ve_0]$ and $r_0>0$ such that 
$
\|\left(F_{\ve,\la}'(u_0)-F'_{\ve,\la}(s u_1 + (1-s) u_2)\right)(u_1 - u_2)\|_U
\le\frac{\al}{2}\|u_1-u_2\|_U
$
for all $\ve \in (0,\ve_1]$, $\la \in \Lambda$, $s \in [0,1]$ and $u_1,u_2 \in \mathbb{B}_{r_0}$.
Hence,
\bee
\label{half}
\|G_{\ve,\la}(u_1) - G_{\ve,\la}(u_2)\|_U
\le \frac{1}{2} \|u_1-u_2\|_U
\mbox{ for all } \ve \in (0,\ve_1],
\la \in \Lambda \mbox{ and } u_1,u_2 \in \mathbb{B}_{r_0}.
\ee

Now, let us show that $G_{\ve,\la}$ maps  $\mathbb{B}_{r_0}$
into  $\mathbb{B}_{r_0}$ for all
$\la \in \Lambda$ and all sufficiently small $\ve>0$.
Take  $\ve\in (0,\ve_1]$, $\la \in \Lambda$ and $u \in  \mathbb{B}_{r_0}$. Then \reff{invest} and \reff{half} imply that
\begin{eqnarray*}
\label{in}
&&\left\|G_{\ve,\la}(u) - u_0\right\|_U \le \left\|G_{\ve,\la}(u) - G_{\ve,\la}(u_0)\right\|_U
+\left\|G_{\ve,\la}(u_0) - u_0\right\|_U\\
&&\le \frac{1}{2}\left\|u-u_0\right\|_U+\left\|(I-F_{\ve,\la}'(u_0))^{-1}(u_0-F_{\ve,\la}(u_0))\right\|_U
\le \frac{r_0}{2}+\frac{1}{\al}\left\|u_0-F_{\ve,\la}(u_0)\right\|_U. 
\end{eqnarray*}
But assumption \reff{limass} yields that, if $\ve_1$ is taken sufficiently small, for all  $\ve \in (0,\ve_1]$ 
and $\la \in \Lambda$
we have 
$\|u_0-F_{\ve,\la}(u_0)\|_U \le \al r_0/2$.
Hence, for those $\ve$ and $\la$ we get $\left\|G_{\ve,\la}(u) - u_0\right\|_U \le r_0$.

Therefore, Banach's fixed point principle yields the following: For all $\ve \in (0,\ve_1]$ and $\la \in \Lambda$
there exists exactly one  $u = u_{\ve,\la} \in
\mathbb{B}_{r_0}$ with $u=F_{\ve,\la}(u)$.

Finally, let us prove \reff{abest}.  We take $\ve \in (0,\ve_1]$ and $\la \in \Lambda$ and estimate as above:
\begin{eqnarray*}
\|u_{\ve,\la}-u_0\|_U &\le& \|G_{\ve,\la}(u_{\ve,\la}) - G_{\ve,\la}(u_0)\|_U+\|G_{\ve,\la}(u_0) - u_0\|_U\\
&\le & \frac{1}{2}\|u_{\ve,\la}-u_0\|_U+\frac{1}{\al}\|u_0-F_{\ve,\la}(u_0)\|_U.
\end{eqnarray*}
Hence, \reff{abest} is true with $\rho=2/\al$.
\qed

\section{Proof of Theorem \ref{main}}
\label{sec3}
\setcounter{equation}{0}
\setcounter{theorem}{0}
In this section we will prove Theorem~\ref{main} by means of Theorem \ref{ift}. Hence, all assumptions of Theorem~\ref{main} (i.e. 
\reff{aass}-\reff{cass}, existence of the weak solution $u=u_0$ to \reff{hombvp}, non-existence of weak solutions $u\not=0$ to \reff{linhombvp})
will be supposed to be satisfied.

Let $r>1$ be fixed.

\begin{lemma}
\label{prep}
For all $\ve>0$ and $B\in \MM$  the following is true:

(i) If $u$ is a weak solution to \reff{ODE}, 
then $u \in W^{1,\infty}((0,1);\R^n)$, and
\bee
\label{ueq}
u(x)=\int_0^x 
\Big(A(\y)+B(\y)\Big)^{-1}\left(\gamma_{\ve,B}(u)-c(y,u(y))
+\int_0^yd(z,u(z))dz\right)dy
\ee
for all $x \in [0,1]$, where
\bee
\label{gammadef}
\gamma_{\ve,B}(u):=\displaystyle\left(\int_0^1
\Big(A(\y)+B(\y)\Big)^{-1}dy\right)^{-1}
\int_0^1\left(c(y,u(y)-\int_0^yd(z,u(z))dz\right)dy.
\ee

(ii) If $u \in C([0,1];\R^n)$ is a solution to \reff{ueq}, then
$u$ is 
a weak solution to~\reff{ODE}. 
\end{lemma} 
{\bf Proof }
Take $\ve>0$ and $B\in \MM$.
 
(i) Let $u\in W_0^{1,2}(\Omega;\R^n)$ be a weak solution to \reff{ODE}. 
Take an arbitrary test function $\vp \in W_0^{1,2}(\Omega;\R^n)$.
Then \reff{vareq} implies that
\begin{eqnarray*}
0&=&\int_0^1\left(\Big(A(\x)+B(\x)\Big)u'(x)+c(x,u(x))\Big)\cdot\vp'(x)-
d(x,u(x))\cdot\int_x^1\vp'(y)dy\right)dx\\
&=&\int_0^1\left(\Big(A(\x)+B(\x)\Big)u'(x)+c(x,u(x))\Big)
-
\int_0^xd(y,u(y))dy\right)\cdot\vp'(x)dx.
\end{eqnarray*}
Therefore $(A(\x)+B(\x))u'(x)+c(x,u(x))\Big)
-
\int_0^xd(y,u(y))dy$ is constant with respect to $x$, i.e.
\begin{eqnarray*}
&&\Big(A(\x)+B(\x)\Big)u'(x)+c(x,u(x))
-\int_0^xd(y,u(y))dy\\
&&=
\int_0^1\left(\Big(A(\y)+B(\y)\Big)u'(y)+c(y,u(y))
-
\int_0^yd(z,u(z))dz\right)dy=:\tilde{\gamma}_{\ve,B}
\end{eqnarray*}
for all $x \in [0,1]$.
It follows that
\bee
\label{ustricheq}
u'(x)=\Big(A(\x)+B(\x)\Big)^{-1}\left(
\tilde{\gamma}_{\ve,B}-c(x,u(x))+\int_0^xd(y,u(y))dy\right)
\ee
for all $x \in [0,1]$,
in particular $u \in W^{1,\infty}(\Omega;\R^n)$. Because of $u(0)=u(1)=0$ it follows also that
$$
0=\int_0^1u'(x)dx=\int_0^1\Big(A(\x)+B(\x)\Big)^{-1}\left(
\tilde{\gamma}_{\ve,B}-c(x,u(x))+\int_0^xd(y,u(y))dy\right)dx,
$$
i.e. $\tilde{\gamma}_{\ve,B}=\gamma_{\ve,B}$
(cf. \reff{gammadef}), and, hence, 
\reff{ustricheq} and the boundary condition $u(0)=0$ imply \reff{ueq}.

(ii) Let $u \in  C([0,1];\R^n)$ be a  solution to \reff{ueq}.
From \reff{ueq} and \reff{gammadef} follows $u(0)=u(1)=0$.
Further, from \reff{ueq} follows that $u\in W^{1,\infty}((0,1);\R^n)$ and
$$
\Big(A(\x)+B(\x)\Big)u'(x)+c(x,u(x))=\gamma_{\ve,B}+\int_x^1d(y,u(y))dy
\mbox{ for a.a. } x \in (0,1).
$$
If we multiply this scalarly  by $\vp'(x)$ with
an arbitrary test function $\vp \in C^1([0,1];\R^n)$ with $\vp(0)=\vp(1)=0$ and integrate with respect to $x$, then we get \reff{vareq}.
\qed\\

Similarly to Lemma \ref{prep} we get
the following: The function $u_0$, which is by assumtion of Theorem \ref{main} a weak solution to the homogenized problem \reff{hombvp},
satisfies
\bee
\label{unulleq}
u_0(x)=A_0^{-1}\int_0^x 
\left(\gamma_0(u)-c(y,u_0(y))
+\int_0^yd(z,u_0(z))dz\right)dy
\ee
for all $x \in [0,1]$, where
\bee
\label{gammanulldef}
\gamma_0(u):=
A_0
\int_0^1\left(c(y,u_0(y))-\int_0^yd(z,u_0(z))dz\right)dy,
\ee
and a function $u$ is a weak solution to the linearized homogenized problem \reff{linhombvp} if and only if
$$
u(x)=A_0^{-1}\int_0^x 
\left(\gamma'_0(u_0)u-\partial_uc(y,u_0(y))u(y))
+\int_0^y\partial_ud(z,u_0(z))
u(z)dz\right)dy
$$
for all $x \in [0,1]$, where
\bee
\label{gammanullstrichdef}
\gamma'_0(u_0)u:=
A_0
\int_0^1\left(\partial_uc(y,u_0(y))u(y)-\int_0^y\partial_ud(z,u_0(z))u(z)
dz\right)dy.
\ee

Now we are going to apply Theorem \ref{ift} in order to solve the boundary value problem \reff{ODE} with $\ve \approx 0$ and $\|u-u_0\|_\infty\approx 0$.
We introduce the setting of Theorem \ref{ift}
as follows:
\bee
\label{abdef}
\left.
\begin{array}{l}
U:=C([0,1];\R^n),\; V:=W^{1,\infty}((0,1);\R^n), \;
\|u\|_U:=\|u\|_\infty,\; \|v\|_V:=\|v\|_\infty+\|v'\|_\infty,\\
\displaystyle[F_0(u)](x):=A_0^{-1}
\int_0^x 
\left(\gamma_0(u)-c(y,u(y))
+\int_0^yd(z,u(z))dz\right)dy,\\
\Lambda:=\MM.
\end{array}
\right\}
\ee
The role of the parameters $\lambda \in \Lambda$ in Theorem \ref{ift} now is played by the localized defects $B \in \MM$, and the role of the maps
$F_{\ve,\la} \in C^1(U;V)$
in Theorem \ref{ift} now is played, for $\ve>0$ and $B \in \MM$,  by
the maps
$F_{\ve,B}:C([0,1];\R^n)\to W^{1,\infty}((0,1);\R^n)$, which are defined by
\bee
\label{abdef1}
[F_{\ve,B}(u)](x):=
\int_0^x 
\Big(A(\y)+B(\y)\Big)^{-1}\left(\gamma_{\ve,B}(u)-c(y,u(y))
+\int_0^yd(z,u(z))dz\right)dy.
\ee
Further, the approximate solution $u_0$
of Theorem \ref{ift} is the solution $u_0$ to the homogenized boundary value problen \reff{hombvp}, which is given by assumption of Theorem \ref{main}.

Because of Lemma \ref{prep} we have the following: 
In order to prove Theorem \ref{main}(i)
we have to verify the assumptions  \reff{compass}-\reff{linlimass}
of Theorem \ref{ift} 
in the setting \reff{abdef}, \reff{abdef1}.
And in order to prove Theorem \ref{main}(ii)
we have additionally to verify the following condition: 
\bee
\label{Feps1}
\left.
\begin{array}{l}
\mbox{If } 
c(\cdot,u_0(\cdot))\in W^{1,\infty}((0,1);\R^n)), 
\mbox{ then }\|F_{\ve,B}(u_0)-F_0(u_0)\|_\infty=O(\ve) \mbox{ for } \ve \to 0\\
\mbox{uniformly with respect to } B \in \MM.
\end{array}
\right\}
\ee

\subsection{Verification of 
\reff{compass}-\reff{linlimass1}}
\label{subsec:Feps}
The Sobolev space $W^{1,\infty}((0,1);\R^n)$ is compactly embedded into $C([0,1];\R^n)$ because of the Arzela-Ascoli theorem. Hence, assumption \reff{compass} of Theorem \ref{ift} is satisfied in the setting \reff{abdef}, \reff{abdef1}.

Assumption \reff{Fnullass}
of Theorem \ref{ift} is satisfied in the setting \reff{abdef}, \reff{abdef1} because of \reff{unulleq} and \reff{gammanulldef}.

The maps $F_{\ve,B}$ and $F_0$ are differentiable from $C([0,1];\R^n)$ into
$W^{1,\infty}((0,1);\R^n)$ because the superposition operators  $u \mapsto c(\cdot,u(\cdot))$ and $u \mapsto d(\cdot,u(\cdot))$ are differentiable
from $C([0,1];\R^n)$ into $L^\infty((0,1);\R^n)$
(cf. assumption \reff{cass}), and
$$
[F'_0(u_0)u](x)
=A_0^{-1}
\int_0^x 
\left(\gamma_0'(u_0)u-\partial_uc(y,u_0(y))u(y)
+\int_0^y\partial_ud(z,u_0(z))
u(z)dz\right)dy
$$
with $\gamma'_0(u_0)u$ defined in \reff{gammanullstrichdef}, and, for $\ve>0$,
\begin{eqnarray*}
&&[F'_{\ve,B}(u_0)u](x)\\
&&=
\int_0^x 
\Big(A(\y)+B(\y)\Big)^{-1}\left(\gamma_{\ve,B}'(u_0)u-\partial_uc(y,u_0(y))u(y)
+\int_0^y\partial_ud(z,u_0(z))
u(z)dz\right)dy
\end{eqnarray*}
with
\begin{eqnarray*}
&&\gamma'_{\ve,B}(u_0)u\\
&&:=\left(\int_0^1
\Big(A(\y)+B(\y)\Big)^{-1}dy\right)^{-1}
\int_0^1\left(\partial_uc(y,u_0(y))u(y)-\int_0^y\partial_ud(z,u(z))
u(z)dz\right)dy.
\end{eqnarray*}

Assumption \reff{linass} of Theorem \ref{ift} is satisfied in the setting \reff{abdef}, \reff{abdef1} because of the assumption of Theorem \reff{main} that there do not exist weak solutions $u \not=0$ to the linearized homogenized problem \reff{linhombvp}.

Further, let us verify assumption \reff{supass} of Theorem \ref{ift} in the setting \reff{abdef}, \reff{abdef1}.
We have to show that
$$
\left\|
\Big(A(\x)+B(\x)\Big)^{-1}\left(\gamma_{\ve,B}'(u_0)u-\partial_uc(x,u_0(x))u(x)
+\int_0^y\partial_ud(y,u_0(y))
u(y)dy\right)\right\|
$$
is bounded uniformly with respect to $\ve>0$, $B \in \MM$, $x \in (0,1)$ and $u \in C([0,1];\R^n)$ with $\|u\|_\infty \le 1$.
But this follows from \reff{cass} and 
\reff{inv}.

Finally, let us verify assumption \reff{linlimass1} of Theorem \ref{ift} in the setting \reff{abdef}, \reff{abdef1}.
We have to show that
\begin{eqnarray*}
&&[(F'_{\ve,B}(u_0+u)-F'_{\ve,B}(u_0))v](x)\\
&&=\int_0^x 
\Big(A(\y)+B(\y)\Big)^{-1}dy
\,\Big(\gamma_{\ve,B}'(u_0+u)v-\gamma_{\ve,B}'(u_0)v\Big)\\
&&\,\;\;\;\;-\int_0^x 
\Big(A(\y)+B(\y)\Big)^{-1}
\Big(\partial_uc(y,u_0(y)+u(y))
-\partial_uc(y,u_0(y))
u(y)dy\\
&&\,\;\;\;\;+\int_0^x 
\Big(A(\y)+B(\y)\Big)^{-1}
\int_0^y\Big(\partial_ud(z,u_0(z)+u(z))
-\partial_ud(z,u_0(z))\Big)
v(z)dzdy
\end{eqnarray*}
tends to zero for $\ve+\|u\|_\infty \to 0$ uniformly with respect to $B \in\MM$, $x \in (0,1)$ and $v \in C([0,1];\R^n)$ with $\|v\|_\infty \le 1$ (in fact it tends to zero for $\|v\|_\infty \to 0$ uniformly with respect $\ve>0$), where
\begin{eqnarray*}
&&\gamma_{\ve,B}'(u_0+u)v-
\gamma_{\ve,B}'(u_0)u\\
&&=
\left(\int_0^1
\Big(A(\y)+B(\y)\Big)^{-1}dy\right)^{-1}
\int_0^1\Big(\partial_uc(y,u_0(y)+u(y))
-\partial_uc(y,u_0(y))\Big)
v(y)dy\\
&&
\,\;\;\;\;-\left(\int_0^1
\Big(A(\y)+B(\y)\Big)^{-1}dy\right)^{-1}
\int_0^y\Big(\partial_ud(z,u(z)+v(z))-
\partial_ud(z,u(z))\Big)
u(z)dzdy.
\end{eqnarray*}
Because of assumption \reff{cass} we have that $u \in\R^n \mapsto  \partial_uc(\cdot,u) \in L^\infty((0,1);\M_n)$ is continuous and, hence, uniformly continuous on bounded sets.
Therefore
$$
\lim_{\|u\|_\infty\to 0}
\sup_{\|v\|\le 1}
\left\|\Big(\partial_uc(\cdot,u_0(x)+u(x))-
\partial_uc(\cdot,u_0(x))\Big)v\right\|_\infty=0
$$
uniformly with respect to $x \in (0,1)$.
Hence
$$
\lim_{\|u\|_\infty\to 0}
\sup_{\|v\|_\infty\le 1}
\left\|\Big(\partial_uc(\cdot,u_0(\cdot)+u(\cdot))-
\partial_uc(\cdot,u_0(\cdot))\Big)v(\cdot)\right\|_\infty=0
$$
And similarly for $\partial_ud$.
Hence, \reff{inv} yields that
$\|\gamma_{\ve,B}'(u_0+u)v-
\gamma_\ve'(u_0)v\|\to 0$ 
for $\|u\|_\infty \to 0$
uniformly with respect to $\ve>0$, $B \in \MM$ and $\|v\|_\infty \le 1$. Hence, again \reff{inv} implies that
$$
\lim_{\|u\|_\infty\to 0}
\sup_{\ve>0,B \in \MM,\|v\|_\infty\le 1}
\|(F'_{\ve,B}(u_0+u)-F'_{\ve,B}(u_0))v\|_\infty=0.
$$

\subsection{Verification of \reff{limass} and
\reff{Feps1}} 
\label{subsec:infass}

The following lemma is the only tool from classical homogenization theory which we are going to use.
For related results see, e.g.
\cite[Proposition 1.1]{Blanc}, 
\cite[Lemma 1.1]{N}, \cite[Proposition 2.2.2]{Shen},
\cite[Lemma 3.1]{Xu}. Roughly speaking, the lemma claims that the homogenized version of the matrix function $A(\cdot/\ve)$ is $A_0$ (cf. \reff{Anulldef}), and 
that the homogenized version of the matrix function $A(\cdot/\ve)+B(\cdot/\ve)$ is $A_0$ also.

\begin{lemma}
\label{prep1}
(i) If $u \in L^1((0,1);\R^n)$, then
\bee
\label{fsup}
\lim_{\ve \to 0}\sup_{0 \le \al \le \beta \le 1, B \in \MM}\left\|
\int_\al^\beta\left(\Big(A(\x)+B(\x)\Big)^{-1}-A_0^{-1}
\right)u(x)dx\right\|=0.
\ee

(ii) There exists $\gamma>0$ such that for all $\ve \in (0,1]$, $B \in \MM$ and $u \in W^{1,\infty}((0,1);\R^n)$ we have
\bee
\label{fsup1}
\sup_{0 \le \al \le \beta \le 1}\left\|
\int_\al^\beta\left(\Big(A(\x)+B(\x)\Big)^{-1}-A_0^{-1}
\right)u(x)dx\right\|\le \gamma \ve \Big(\|u\|_\infty+\|u'\|_\infty\Big).
\ee
\end{lemma}
{\bf Proof }
(i) We proceed as in \cite[Proposition 1.1]{Blanc}. 
Because of \reff{inv} the map
$$
u \in L^1((0,1);\R^n)\mapsto \sup_{1 \le \al \le \beta \le 1}\left\|
\int_\al^\beta\left(\Big(A(\x)+B(\x)\Big)^{-1}-A_0^{-1}
\right)u(x)dx\right\|\in \R
$$
is continuous uniformly with respect to $\ve>0$ and $B \in \MM$. Moreover, the set of all piecewise constant functions is dense in $L^1((0,1);\R^n)$.
Therefore it sufficies to prove \reff{fsup} for piecewise constant functions $u$, i.e. to prove
$$
\lim_{\ve \to 0}\sup_{0 \le \al \le \beta \le 1, B \in \MM}\left\|
\int_\al^\beta\left(\Big(A(\x)+B(\x)\Big)^{-1}-A_0^{-1}
\right)dx\; v\right\|=0
\mbox{ for all } v \in \R^n.
$$
Hence, it remains to prove assertion (ii) of the lemma.

(ii) Take $\ve>0$, $B \in \MM$ and $u \in W^{1,\infty}((0,1);\R^n)$. Because of \reff{inv} and of $\|B\|_1 \le r$ (cf. \reff{Mdef}) we have that
\begin{eqnarray*}
&&\left\|\int_\al^\beta\left(\Big(A(\x)+B(\x)\Big)^{-1}-A(\x)^{-1}
\right)u(x)dx\right\|\nonumber\\
&&=\left\|\int_\al^\beta\left(\Big(A(\x)+B(\x)\Big)^{-1}B(\x)A(\x)^{-1}
\right)u(x)dx\right\|\nonumber\\
&&=\ve\left\|\int_{\al/\ve}^{\beta/\ve}\left(\Big(A(y)+B(y)\Big)^{-1}B(y)A(y)^{-1}
\right)u(\ve y)dy\right\|\le \mbox{const }\ve
\|u\|_\infty,
\end{eqnarray*}
where the constant does not depend on $\ve$, $B$, $\al$, $\beta$ and $u$.
Hence, in order to prove \reff{fsup1} it remains to prove that
\bee
\label{fsup4}
\sup_{0 \le \al \le \beta \le 1}\sup_{\|u\|_\infty+\|u'\|_\infty\le 1}\left\|
\int_\al^\beta\Big(A(\x)^{-1}-A_0^{-1}
\Big)u(x)dx\right\|=O(\ve)
\mbox{ for } \ve \to 0.
\ee

Define $\A \in L^\infty(\R;\M_n))$ by
$
\A(y):=A(y)^{-1}-A_0^{-1}.
$
Then 
$\A(y+1)=\A(y)$ and, because of \reff{Anulldef}, $\int_y^{y+1}\A(z)dz=0$ for all
$y \in \R$, and \reff{inv} yields that
$$
\gamma_0:=\mbox{ess sup}\{\|\A(y)v\|:\; y \in \R,\; v \in \R^n;\; \|v\|\le 1\}<\infty.
$$
Take $1 \le \al < \beta \le 1$, $u \in \R^n$ with $\|u\|\le 1$, and take $\ve>0$ sufficiently small.
Then
\begin{eqnarray*}
&&\left\|\int_{\al}^{\beta}\left(A(\x)^{-1}-A_0^{-1}
\right)u(x)dx\right\|=\ve \left\|\int_{\al/\ve}^{\beta/\ve}\A(y)u(\ve y)dy\right\|\\
&&=\ve\left\|\sum_{j=[\al/\ve]+1}^{[\beta/\ve]-1}\int_{j}^{j+1}\A(y)
\Big(u(\ve y)-u(\ve j)\Big)dy+
\int_{\al/\ve}^{[\al/\ve]+1}
\A(y)u(\ve y)dy +
\int_{[\beta/\ve]}^{\beta/\ve}
\A(y)u(\ve y)dy\right\|\\
&&\le \gamma_0 \ve
\left(\sum_{j=[\al/\ve]+1}^{[\beta/\ve]-1}\int_{j}^{j+1}
\|u(\ve y)-u(\ve j)\|dy+2\|u\|_\infty\right).
\end{eqnarray*}
Here $[\al/\ve]$ and $[\beta/\ve]$ are the integer parts of $\al/\ve$ and $\beta/\ve$, respectively.
But
$$
\int_{j}^{j+1}
\|u(\ve y)-u(\ve j)\|dy\le\ve\;
\|u'\|_\infty\int_{j}^{j+1}(y-j)
dy=\frac{\ve}{2}\;\|u'\|_\infty,
$$
therefore
$$
\sum_{j=[\al/\ve]+1}^{[\beta/\ve]-1}\int_{j}^{j+1}
\|u(\ve y)-u(\ve j)\|dy
\le \Big([\beta/\ve]-[\al/\ve]-2\Big)\;
\frac{\ve}{2}\;\|u'\|_\infty
\le \mbox{const } \|u'\|_\infty,
$$
where the constant does not depend on $\ve$, $\al$, $\beta$ and $u$.
Hence, \reff{fsup4} is proved.
\qed\\

\begin{remark}
\label{Buniform}
The proof of Lemma \ref{prep1} shows that the constant $\gamma$ in \reff{fsup1} depends on the matrix functions $A$ and $B$ via the norms of $\|A\|_\infty$, $\|B\|_\infty$ and  
$\|B\|_1$ and via the suprema in \reff{inv}, only.
\end{remark}

\begin{remark}
\label{unifo}
For $\ve>0$ and $B \in \MM$ define
$$
M_{\ve,B}:=\left(\int_0^1
\Big(A(\x)+B(\x)\Big)^{-1}dx\right)^{-1}.
$$
Then Lemma \ref{prep1}(ii) yields that
$\left\|M_{\ve,B}^{-1}-A_0^{-1}\right\|=O(\ve)$ for  $\ve \to 0$ uniformly with respect to $B \in \MM$.
Hence, also
$\left\|M_{\ve,B}-A_0\right\|=O(\ve)$ for  $\ve \to 0$ uniformly with respect to $B \in \MM$.
Because of 
$$
\gamma_{\ve,B}(u_0)-\gamma_0(u_0)=\left(M_{\ve,B}
-A_0\right)
\int_0^1\left(c(x,u_0(x))-\int_0^xd(y,u_0(y))dy\right)dx
$$
and, for $u \in C([0,1];\R^n)$, of
$$
\gamma'_{\ve,B}(u_0)u-\gamma'_0(u_0)u
=\left(M_{\ve,B}
-A_0\right)
\int_0^1\left(\partial_uc(x,u_0(x))u(x)-\int_0^x\partial_ud(y,u_0(y))u(y)dy
\right)dx
$$
it follows that
\bee
\label{gammalim}
\left.
\begin{array}{l}
\|\gamma_{\ve,B}(u_0)-\gamma_0(u_0)\|
+\sup_{\|u\|_\infty \le 1}\|\gamma_{\ve,B}'(u_0)u-\gamma'_0(u_0)u\|=O(\ve)
 \mbox{ for } \ve \to 0\\
 \mbox{uniformly with respect to } B \in \MM.
 \end{array}
 \right\}
\ee
\end{remark}

Now, let us verify \reff{limass}
of Theorem \ref{ift} 
in the setting \reff{abdef}, \reff{abdef1}. 
Take $\ve>0$ and $B \in \MM$.
Because of the definitions of the maps $F_{\ve,B}$ and $F_0$ in \reff{abdef} we have for all $x \in [0,1]$ that
\begin{eqnarray}
&&[F_{\ve,B}(u_0)](x)-[F_0(u_0)](x)\nonumber\\
&&=\int_0^x 
\Big(A(\y)+B(\y)\Big)^{-1}\left(\gamma_{\ve,B}(u_0)-c(y,u_0(y))
+\int_0^yd(z,u_0(z))dz\right)
dy\nonumber\\
&&\;\;\;\;-
A_0^{-1}
\int_0^x 
\left(\gamma_0(u_0)-c(y,u_0(y))
+\int_0^yd(z,u_0(z))dz\right)
dy\nonumber\\
\label{Fest}
&&=\int_0^x 
\Big(\Big(A(\y)+B(\y)\Big)^{-1}-A_0^{-1}\Big)
\left(\gamma_{\ve,B}(u_0)-c(y,u_0(y))
+\int_0^yd(z,u_0(z))dz
\right)dy\nonumber\\
&&\;\;\;\;+x A_0^{-1}\Big(\gamma_{\ve,B}(u_0)-\gamma_0(u_0)\Big).
\end{eqnarray}
For $y \in \R$ denote
\bee
\label{Avedef}
M_{\ve,B}(y):=A(\y)+B(\y).
\ee
Then \reff{gammalim} and Lemma \ref{prep1}(i) yield that
\begin{eqnarray*}
&&\lim_{\ve\to 0}\|F_{\ve,B}(u_0)-F_0(u_0)\|_\infty\\
&&=\lim_{\ve\to 0}\sup_{x \in [0,1]}
\left\|\int_0^x 
\Big(M_{\ve,B}(y)^{-1}-A_0^{-1}\Big)
\left(\gamma_0(u_0)-c(y,u_0(y))
+\int_0^yd(z,u_0(z))
dz\right)dy\right\|=0
\end{eqnarray*}
uniformly with respect to $B \in \MM$.
Hence,  assumption \reff{limass}
of Theorem \ref{ift} 
in the setting \reff{abdef}, \reff{abdef1}
is verified.

And \reff{gammalim}, \reff{Fest} and Lemma \ref{prep1}(ii) yield that condition
\reff{Feps1} is verified also.

\subsection{Verification of \reff{linlimass}} 
\label{subsec:infass1}
Take $\ve>0$, $B \in \MM$ and $u \in C([0,1];\R^n)$. Using notation 
\reff{Avedef} again, we get for all $x \in [0,1]$ that
\begin{eqnarray*}
&&[F'_{\ve,B}(u_0)u](x)-[F'_0(u_0)u](x)\\
&&=\int_0^x 
\Big(M_{\ve,B}(y)^{-1}-A_0^{-1}\Big)
\left(\gamma'_{\ve,B}(u_0)u-\partial_uc(y,u_0(y))u(y)
+\int_0^y\partial_ud(z,u_0(z))u(z)dz
\right)dy\nonumber\\
&&\;\;\;\;+x A_0^{-1}\Big(\gamma'_{\ve,B}(u_0)u-\gamma'_0(u_0)u
\Big).
\end{eqnarray*}
Hence, \reff{gammalim} and Lemma \ref{prep1}(i) yield
\begin{eqnarray*}
&&\lim_{\ve\to 0}\|(F'_{\ve,B}(u_0)-F'_0(u_0))u\|_\infty\\
&&=\lim_{\ve\to 0}\sup_{x \in [0,1]}
\left\|\int_0^x 
\Big(M_{\ve,B}(y)^{-1}-A_0^{-1}\Big)
\left(\partial_uc(y,u_0(y))u(y)
-\int_0^y\partial_ud(z,u_0(z))u(z)
dz\right)dy\right\|=0
\end{eqnarray*}
uniformly with respect to $B \in \MM$.

\begin{remark}
\label{uniform}
The last limit above is not 
uniform with respect to $u \in C([0,1];\R^n)$ with $\|u\|_\infty \le 1$, i.e. the linear operators $F'_{\ve,B}(u_0)$ tend to $F'_0(u_0)$
for $\ve \to 0$ strongly, but not with respct to the uniform operator norm in ${\cal L}(C([0,1];\R^n))$. This is also the case if the matrix functions $A$ and $B$ and the vector functions $c$ and $d$ are smooth.
\end{remark}

\section*{Acknowledgments}

No founds or grands where received for the submitted work. The author has no relevant financial
or non-financial interests to disclose.
There is no data used.


\begin{thebibliography}{troi}

\bibitem{Ben}
A. Bensoussan, J.L. Lions, G. Papanicolaou, Asymptotic Analysis for Periodic Structures. Studies in Mathematics and its  Applications vol. {\bf 3}, North-Holland, 1978. 

\bibitem{Blanc}
X. Blanc, C. Le Bris,
Homogenization Theory for Multiscale Problems. An Introduction. 
Modeling, Simulation and Applications vol. {\bf 21}, Springer, 2023.




\bibitem{Bun} R. Bunoiu, R. Precup, Localization and multiplicity in the homogenization of nonlinear problems.
Adv. Nonlinear Anal. {\bf 9} (2020), 292--304.




\bibitem{Butetc} V.F. Butuzov, N.N. Nefedov, O.E. Omel'chenko, L. Recke,
Time-periodic boundary layer solutions to singularly perturbed parabolic problems.
J. Differ. Equations {\bf 262} (2017), 4823--4862.

\bibitem{But2022}  
V.F. Butuzov, N.N. Nefedov,  O.E. Omel'chenko, L. Recke, 
Boundary layer solutions to singularly perturbed quasilinear systems.
Discrete Cont. Dyn. Syst., Series B {\bf 27}
(2022), 4255--4283.


\bibitem{Fiedler}  
V.F. Butuzov, N.N. Nefedov,  O.E. Omel'chenko, L. Recke, K.R. Schneider, 
An implicit function theorem and applications to nonsmooth boundary layers.
In: {\it Patterns of Dynamics}, ed. by 
P. Gurevich, J. Hell, B. Sandstede, A. Scheel, Springer Proc. in Mathematics \& Statistics vol. {\bf 205}, Springer,  
2017, 111--127.
















\bibitem{He}  
Wen Ming He, Jun Zhi Cui, Error estimate of the homogenization solution for elliptic problems with small periodic coefficients on $L^\infty(\Omega)$. Science China Mathematics {\bf 53} (2010), 1231--1252.




\bibitem{Katz}
S.G. Krantz, H.R. Parks, The Implicit Function Theorem. History, Theory, and Applications.
Birkh\"auser, 2002.




\bibitem{Lanza1} M. Lanza de Cristoforis, P. Musolino, Two-parameter homogenization for a nonlinear periodic Robin problem for a Poisson equation: a functional analytic approach. Rev. Mat. Complut. {\bf 31} (2018), 63--110.

\bibitem{Lanza2} M. Lanza de Cristoforis, P. Musolino, Asymptotic behaviour of the energy integral of a two-parameter homogenization
problem with nonlinear periodic Robin boundary conditions. Proc. Edinb. Math. Soc. II. Ser. {\bf 62} (2019), 985--1016.




\bibitem{Magnus1} R.J. Magnus, The implicit function theorem and multi-bump solutions of periodic partial differential equations. Proc.
Royal Soc.  Edinb.  {\bf 136A} (2006), 559--583.


\bibitem{Magnus2} R.J. Magnus, A scaling approach to bumps and multi-bumps for nonlinear partial differential 
equations. Proc.
Royal Soc.  Edinb.  {\bf 136A} (2006), 585--614.


\bibitem{NURS}  
N.N. Nefedov, A.O. Orlov, L. Recke, K.R. Schneider, Nonsmooth regular perturbations of singularly perturbed problems. 
J. Differ. Equations {\bf 375} (2023), 206--236.


\bibitem{I} N.N. Nefedov, L.~Recke, A common approach to singular perturbation and homogenization I: Quasilinear ODE systems. arXiv:2309.15611.


\bibitem{II} N.N. Nefedov, L.~Recke, A common approach to singular perturbation and homogenization II: Semilinear elliptic PDE systems. J. Math. Anal. Appl. {\bf 545} (2025), Article ID 129099.


\bibitem{N}  
S. Neukamm, An introduction to the qualitative and quantitative theory of homogenization. Interdisciplinary Information Sciences {\bf 22} (2016), 147--186.



\bibitem{OmelchenkoRecke2015} O.E.~Omel'chenko, L.~Recke,
Existence, local uniqueness and asymptotic approximation of spike solutions
to singularly perturbed elliptic problems. Hiroshima Math. J. {\bf 45} (2015), 35--89.



\bibitem{Recke2022} L.~Recke, Use of very weak approximate
boundary layer solutions to spatially nonsmooth
singularly perturbed problems.
J. Math. Anal. Appl. {\bf 506} (2022), Article ID 125552.

\bibitem{ReckeOmelchenko2008} L.~Recke, O.E.~Omel'chenko,
Boundary layer solutions to problems
with infinite dimensional singular and regular perturbations.
J. Differ. Equations {\bf 245} (2008), 3806--3822.

\bibitem{Shen}
Zongwei Shen, Periodic Homogenization of Elliptic Systems. Operator Theory: Advances and Applications vol. {\bf 269}, Birkh\"auser, 2018.







\bibitem{Xu} Shixin Xu, Xingye Yue, Changrong Zhang,
Homogenization: in mathematics or physics? Discrete Contin. Dyn. Syst. Series S {\bf 9} (2016), 1575--1590.

\bibitem{YXu} Yao Xu, Weisheng Niu, Homogenization of elliptic systems with stratified structure revisited. Commun. Partial Differ. Equations
{\bf 45} (2020), 655--689.





\end{thebibliography}
\end{document}